\documentclass[12pt]{article}

\usepackage{graphicx}
\usepackage{xcolor}
\usepackage{float}
\usepackage{url}
\usepackage{natbib}
\usepackage[normalem]{ulem}

\usepackage{epstopdf}
\usepackage{hyperref} 
\hypersetup{colorlinks=true, linkcolor=blue,  anchorcolor=blue,  
citecolor=blue, filecolor=blue, menucolor=blue, pagecolor=blue,  
urlcolor=blue,pdftitle={ACmain}}

\title{ \textbf{Community, Collaboration, and Climate}}
\author{Johanna Hardin and Shahriar Shahriari}
\date{April 2022}

\begin{document}

\maketitle

\abstract{The Department of Mathematics \& Statistics at Pomona College has long worked to create an inclusive and welcoming space for all individuals to study mathematics.  Many years ago, our approach to the lack of diversity we saw in our majors was remediation through programming which sought to ameliorate student deficits.  More recently, however, we have taken an anti-deficit approach with focus on changes to the department itself.  The programs we have implemented are described below as enhancing community, collaboration, and climate within our department.}\\

\noindent {\bf Keywords:} anti-deficit; diversity, equity, and inclusion; access;

\section{Introduction}

The murders of George Floyd,  Breonna Taylor, Ahmaud Arbery, Tony McDade, Philando Castile, Tamir Rice, Michael Brown, Eric Garner, Trayvon Martin, and countless others, as well as the drastic disparities in death rates due to COVID-19 between whites and communities of color has brought to sharp focus the systematic and structural racism in the US. As a related byproduct, the importance of  providing access to a quality education for all is ever more self-evident. Our math and stats department---like many others---has struggled with its role in providing an inclusive and inviting climate for students from groups that are traditionally underrepresented in mathematics. With the hope that small and concrete changes can make a difference, and that sharing our experiences may be helpful to others, we are eager to join in the discussion. While we still have a long way to go, we have, over the years, implemented a number of concrete programs, that have been positive overall. 
We are not speaking as the voice of our department, but we hope that our personal perspective on the modest and concrete programs in our department over the past decade can be of some use to others.

Consider the following two different perspectives on the job of a mathematics or statistics faculty at a selective undergraduate college: (1) Our main job is to identify talent, to encourage the talented to pursue the mathematical sciences, and to lead the rest to other, more ``appropriate", disciplines, and (2) we should have the philosophy that \emph{everyone} can do math, and that our job is to facilitate \emph{every} student's  success in mathematics. In the math and stats department at Pomona College, we never explicitly discussed these points of view, but some variation of the tension between them has always been in the background, and, if you take a bird's eye view of the last thirty years, despite the range of personalities and continued disagreements, we have slowly but decisively moved from the first point of view toward the second. 

The opening dichotomy is, of course, somewhat rhetorical. Ideally, students find success in mathematics at different levels; some continue in the mathematical sciences and some don't. What is not rhetorical is the data that, over a long period of time, we had very little success in attracting students of color (particularly Black, Hispanic, and indigenous students), low-income students, first generation students, and even women.  

Pomona College is a small liberal arts college with enviable resources, our faculty do their best to be active researchers while being dedicated teachers and mentors, and our student to faculty ratio is low. Even so, consistently and for many years, the department had attracted fewer students of color into the major than their proportional representation at the College. Even more alarming was the fact that often students of color, students from low income backgrounds, and first-generation students did poorly in our beginning courses. It is not surprising that such students, by and large, did not continue taking mathematics, and steered clear of STEM fields. Whether we liked it or not, we acted as an effective filter across STEM.

Our first response, starting about three decades ago, was curricular. We worried about the students' prior preparation, the quality of their high school work, and the mismatch with our courses. Various faculty developed new courses that relied less on students' prior courses, were interesting, strengthened problem solving skills, and prepared students for future work.\footnote{Our colleague Erica Flapan's \emph{Problem Solving for the Sciences}---a course organized around doing non-calculus ``word problems''---and the second author's \emph{Approximately Calculus}---later turned into a book---were attempts in this direction.}
Many of us introduced or more fully embraced active student learning---with or without technology---in our classes. We reorganized the curriculum to create multiple pathways through the introductory material.\footnote{Thirty years ago, at Pomona, if you wanted to do math, you took three semesters of calculus followed by two semesters of a combined linear algebra/differential equations class before you had any choice of topics or classes.} All of these changes were, on balance, positive. The students responded, and our major and courses grew. The problem of too few students of color, low-income students, and first generation students who continued on to the next math class, however, persisted.

About ten years ago, two students of color, who had done superbly in the second author's own honors calculus class, were failing his linear algebra class. This was one among many data points suggesting that a focus on student preparation cannot possibly be enough. It seemed to us that background preparation was not a determining factor as students of color with all different backgrounds were underperforming what their background preparation alone would have predicted.  The college's work with the Posse Foundation and the work of Uri Treisman\footnote{Our colleague Ami Radunskaya had worked directly with Uri Treisman, the second author was Pomona's inaugural Posse faculty mentor, and the first author has advised a Posse cohort also.} \citep{asera2001} helped us shift our focus to the ``climate'' in the department. We moved the spotlight from the students and their preparation to us, the faculty, and the environment in the department. We slowly started to understand that some of our students did not feel that they belonged, and, as a result, they were not using the available resources (mentor sessions, office hours, tutors, etc.) as effectively as possible. The stress of tackling a difficult discipline in an unfriendly environment---which extended far beyond the math and stats department---adversely affected the students' successes in the classroom and their desire to continue in mathematics. The conclusion was that we had to create extracurricular programs to challenge the dominant climate in the department.
 
At this point it is worth clarifying that while our goal is to reach students who are minoritized or marginalized for any of a variety of reasons, generally, depending on the program, we include students across the board who self-identify as an individual who wants to benefit from the programming.  That is, our programs are populated with students who do not see themselves in mathematics (are minoritized); feel dismissed in the mathematical space (are marginalized); or want to learn and grow in a supportive environment designed to create equity in mathematics.  With our approach, the vast majority (but not all) of students in our programs are students of color, low-income students, and/or first generation students.  Is our approach a strength or a weakness?  We don't have the answer.  But we do know that we firmly believe that mathematics is for everyone, and our inclusive attitude is integrated in to all that we do.

\section{What have others done?}

As posed above, there are generally two approaches to broadening participation in the mathematical sciences.  Substantial work has been done to implement creative ideas in the college mathematics classroom surrounding each of these.  Our paper is not meant to be a comprehensive review of different programming, but we provide a few examples below to contextualize our work with that of the larger mathematics community.

\subsection{A deficit approach}

A deficit approach is typically connected to the amount of previous {\em mathematical training} a student has had prior to their college mathematics classes.  \citet{peck2021} provides a comprehensive overview of the problems surrounding a deficit model, including a succinct definition of a deficit perspective:

\begin{quote}
{\bf Definition:} The deficit perspective is a propensity to locate the source of academic problems in deficiencies within students, their families, their communities, or their membership in social categories (such as race and gender).
\end{quote}

One possible reaction to the definition is that it indicates a personal failing of the student which teachers of mathematics surely do not mean to imply.  However, from the perspective of a student of mathematics, an insufficient background {\em is} felt as a personal failing.  

Online programs have potential to fill in gaps of understanding either throughout the semester (remediation through online just-in-time teaching in order to prevent pre-calculus and calculus from being barriers toward all the STEM fields \citep{bertrand2021}) or prior to college in order to avoid taking remedial mathematics classes (using the online software ALEKS during high school to develop mathematical skills  \citep{fine2009}).

In-person programs can also be effective.  A one-credit supplemental co-requisite Calculus I course \citep{hancock2021} allows for active learning in the classroom as well as reinforcement of the calculus concepts and community within the course.  As has been well-documented, working with faculty on research projects is associated with retention \citep{hernandez2018} and achievement \citep{parker2018}.

\subsection{An anti-deficit approach}

An anti-deficit approach tackles issues which are outside of the individual student's purview.  Again, as defined by  \citet{peck2021}:

\begin{quote}
{\bf Definition:} Anti-deficit perspectives locate the source of academic problems within institutional structures
that work to limit access to educational opportunities. The focus is on the assets that students bring to the
classroom, rather than what they lack.
\end{quote}

\citet{peck2021} provides an extensive review of issues and practical suggestions toward creating anti-deficit mathematics programs.  Additionally, \citet{hagman2021} identifies diversity, equity, and inclusion as one of the core components which make up a successful calculus program.  Working through a classroom activity, an anti-deficit approach to teaching is developed by \citet{adiredja2019}.

In retrospect, many of the steps we have taken at Pomona College to create a math and stats department that is more diverse, equitable, and inclusive follow an anti-deficit model.  We describe these next.

\section{What did we implement?}

The Department of Mathematics \& Statistics has always valued community. Since the early 1990s, we have organized, for all of our lower division classes, evening peer mentoring sessions where students work together on their assignments under the supervision of a more experienced undergraduate. These mentor sessions have encouraged collaboration and have resulted in a buzz of mathematical activity in our building in the evenings. As a department, we encourage collaboration and want the students to compete with the material instead of with each other. Many of our faculty assign challenging projects that are meant for collaborative work, and mentor sessions play a crucial role in bringing the students together to collaborate. But student experience at the mentor sessions has been very uneven. While many of our students loved collaborating in the mentor sessions, and credited these evening sessions with a strong sense of community in the department, others felt left out. Somewhat unsurprisingly, students of color, low-income students, and first generation students were more often than not in the latter group. They felt out of place, did not participate in the discussions, were not comfortable asking questions, and often found themselves isolated.

Here, we report on three on-going programs---two in the department and one college wide---that we believe have had a positive impact on both our department and the students we are hoping to target. The first one, \emph{Learning Communities} (abbreviated LiCMath) are organized for individual classes. The second one is a department-wide cohort program called the \emph{Pomona Scholars of Mathematics} (PSM), and the third one, \emph{Pomona Academy for Youth Success} (PAYS) is a college wide outreach program for local high school students. Our experience with these programs are of varying lengths. The first group of PAYS students arrived on campus in the summer of 2003. The first LiCMath group was organized in the fall of 2010, and the first cohort of PSM students began in the fall of 2014.  We are not presenting an empirical study (as the programs were not set up with a particular experimental design); instead, we present the programs and the evidence we have that they are working to bring about change in our department.

\subsection{LiCMath} 

The literature describes learning community structures to support mathematics programs.  However, while aligned in goals of improved mathematics content delivery, the programs often have different goals in terms of who and how they serve the students.  For example, learning communities play a role in encouraging women in mathematics as described by \citet{evans2018}.  Linking mathematics and English courses, \citet{piercey2017} have found learning communities to be effective at improving problem solving; linking biology and statistics courses \citet{ryu2022} have found that learning communities improve learning across both disciplines.   General college sponsored learning communities can be quite effective {\bf and} efficient when using peer mentors \citep{rieske2015}.  A set of qualitative interviews demonstrated that (college wide)  learning communities help students interact at many different levels and that the first year was an ideal time for learning community involvement \citep{firmin2013}. 

At Pomona College in the Department of Mathematics \& Statistics, a learning community is a small group of students who are all enrolled in one particular class.  They commit to work together and to support each other throughout the semester. A learning community is meant to become a supportive and safe space to work on homework problems, to discuss the lectures, and to bounce ideas off each other.  The meetings take place outside of class times, and, at least initially, under the supervision of an upper class undergraduate student mentor. 

The targeted students are those who do not feel at home in the math and stats department---or in the college in general---and have doubts about whether they belong. Such students could have a strong mathematics background or a weak one but may try to go through the class in isolation. They may, consciously or unconsciously, believe that using the available resources (faculty office hours, mentor sessions, discussions with peers) is a sign of weakness and that they should be able to succeed by sitting alone in their room and working doubly hard on the mathematics. It may not be easy for such students to create a safe collaborative learning environment early in the semester. Such students can come from all walks of life, but first-generation college students, women, students of color, students from less prestigious high schools, and students from less privileged families tend to fall in this category most often. Learning communities are open to all students, but we make a special effort to encourage the students in our targeted audience to participate. 

A faculty member who decides to organize a learning community in their class, announces the aim and structure of LiCMath to the whole class and invites everyone who is interested to apply. In addition, targeted students are sent customized individual emails encouraging them to participate. In our pitch, we consistently emphasize the value of collaboration, and we explain that every student needs a supportive group of peers in order to get the most out of a mathematics class. Some students already have such a support group; others will find it in their dorms or in the mentor sessions. LiCMath is an alternative for students who want a more structured collaborative environment. We explicitly make it clear that students can participate in LiCMath regardless of their perception of their strengths or weaknesses in mathematics. All applicants who apply and who accept the responsibilities of membership are accepted and organized into small groups with 3 to 7 members each. 

The learning communities meet for at least two hours a week, and the students commit to going to all of the meetings regardless of whether they have specific course related questions. If a student believes that they already know the material for that week, then they are still expected to attend the weekly meeting to provide support for others. To build community and to start getting to know each other, we suggest that the students begin each meeting with a short check-in about their week. We also encourage the students not to focus solely on the next assignments. Rather, they should spend some time going over the main points of the material recently covered in class, and discuss their returned homework assignments. By communicating to each other what is interesting, important, challenging, or unclear, they start the process of learning from each other.

While, on paper, it seems as though LiCMath adds hours of work each week to a student's plate, we almost never get feedback that LiCMath is a burden.  Instead, the time spent in a learning community typically replaces time that a student would have been working alone, possibly frustrated about not making progress.  Additionally, working in the learning community increases both confidence and communication skills, making the class mentor sessions much more efficient for the student.  

Our concept of a learning community owes much to the Posse Foundation's idea of creating a peer support group for student success. It also builds on Uri Treisman's pioneering ideas for ``Emerging Scholars'' workshops for calculus classes. We are a small liberal arts college where most students know each other, and where most of our classes are not overburdened with students from other departments simply fulfilling requirements, and, hence, already incorporate many challenging projects and experiences. Our specific circumstances have meant some possibly important modifications to previous models. For example, our learning communities are not marketed as ``honors'' sections, and they do not usually incorporate additional more challenging course material. Nor are the learning communities remedial tutoring sessions for students who are falling behind in the class.  Our students are told that they should join learning communities because isolation is not a path to success in mathematics, and everyone needs a supportive group of peers through the semester to be able to get the most out of any of our classes.

Learning communities are a relatively low-cost structure for creating community. The interested faculty member initiates the process, recruits the students and the mentor, and provides some early guidance. While it takes some time to think through the organization, and striking the right tone in communicating with the students can sometimes be challenging, there is very little on-going work for the faculty during the semester.

We count our students' participation and even enthusiasm for LiCMath as circumstantial evidence for its success. Initially, LiCMath sessions were organized only for lower division classes, but, in a number of upper division classes, the students themselves asked for and organized LiCMath as well. In fact, both Pomona  Physics and Chemistry departments started experimenting with learning communities because students, who had been part of LiCMath, asked for them.

\subsection{PSM} The first semester at college, especially as a full time student living away from family, can be very stressful. Learning to navigate the expectations and the resources of a university while adjusting to a new social landscape is challenging. This is especially so for first generation college students or for individuals who find the physical space, the bureaucratic structure, the social scene, or the academic expectations unfamiliar. To counter the sense of isolation in facing the new environment, and to level the academic playing field, Pomona College has been experimenting with creating cohorts of incoming students centered around academic disciplines. As a part of our efforts, the authors, beginning in the fall of 2014, started the Pomona Scholars in Mathematics (PSM). The goal of PSM is two fold. On the one hand, we want to create a community of students who experience the first year of college and its challenges and opportunities together, and, on the other hand, we hoped that, over time, there would be a spill over effect, and the sense of inclusion and community in the whole department would be enhanced.   

Recognizing the special role of mathematics as a gatekeeper or conduit to other fields, PSM aims to support students who express an interest in fields where mathematics plays an important role (math, statistics, computer science, physics, economics, engineering, etc). In the summer before their first year at Pomona, targeted students---students of color, low-income students, and first generation students---are invited to join PSM. Approximately 15 students do so, and they are advised by two faculty members (usually from the math and stats department but sometimes from aligned disciplines like computer science or physics). The program is built on three pillars:

\begin{itemize}
\item\ {\bf early intensive one-on-one faculty-student advising}.  Each of the two faculty members are assigned to be the faculty academic advisor for half of the students, and the faculty have biweekly one-on-one meetings with their PSM advisees throughout the first year. Taking a wholistic approach to mentoring, these meetings go beyond just academics. We check in with the students about their college experience and home life as well as how their classes are going.

\item\ {\bf weekly group lunch meetings}. The PSM group of first year students and faculty mentors meet once a week over lunch for the duration of the students' first year in college. While the primary purpose of these meetings is to build community and camaraderie---we routinely have the students check in and share their ongoing challenges, successes, and insights---we also introduce the students to college resources, and we regularly invite faculty members from a variety of disciplines to come in and share their academic journey. Our hope is that demystifying the academy and normalizing the students' experiences help counter the sense of isolation in a new and challenging academic environment.

\item\ {\bf within-course clusters}. Our PSM students have a wide range of mathematical background, and they do not all start in the same first mathematics course as first year students. Even so, we encourage the PSM students to take classes together when they can, and, in keeping with our mathematics-is-a-collaborative-enterprise mantra, we urge them to form small study clusters in various sections of math. 
\end{itemize}
After the first year, the PSM students join what we call PSM+. PSM+ is a loosely organized community of sophomore, junior, and senior students that aims to continue the goals of PSM. The larger group meets about once a month over lunch, and organizes other department-wide events.

While it is difficult to tease out confounding variables, our impression is that PSM has been successful in changing the feel of the department. Not only are more of the targeted students continuing in mathematics courses, but more students of color, low-income students, and first generation students have been taking on leadership roles in the department, and, in general, have been more visible as part of Pomona's mathematics community.

Unlike learning communities, the downside of PSM is the need for heavy faculty commitment. Weekly group meetings and biweekly advising meetings are extremely time consuming. The faculty members receive a small stipend for their work but no other compensation. 

Of all our programs, PSM is the one which embraces an anti-deficit approach most fully.  In our weekly meetings we center the experience of the PSM students and who they are as a whole person.  In community, they are able to see that their experiences are sometimes not unique, and they are able to process how to be successful in college.  We also bring in older students, alumni, and other role models (e.g., staff and faculty) who have had non-traditional paths.  Although the majority of guests (and students) are POC, their more important connection is that to a person, they've had to overcome systemic barriers.  Our goal in centering their experiences is to have them see themselves being successful, in their {\bf own} way.   The PSM community empowers our students to understand and challenge the oppressive structures they encounter (see Section \ref{kayla}).

\bigskip

\subsection{PAYS} Pomona Academy for Youth Success (PAYS\footnote{https://www.pomona.edu/administration/draper-center/pays}) is a Pomona College summer academic program for high school students. The program, which began in 2003, has an enviable record of its participants continuing onto selective colleges and universities. The program brings 30 rising tenth graders from local under-served high schools to spend 4 weeks in residence at the College. These same students come back for another 4 weeks as rising eleventh graders and again as rising twelfth graders. 

Each four week session is intensely academic, and its core consists of each cohort taking two classes that are taught by regular college faculty: a math class and a humanities/social sciences class. The classes meet four days a week in the mornings. Afternoons are focused on working collaboratively on homework, workshops on college application processes, and enrichment activities. The faculty are aided by a group of dedicated undergraduate students who serve as teaching assistants to the faculty and as mentors to the students.

 The Department of Mathematics \& Statistics has been instrumental in the design and the implementation of the PAYS program.\footnote{The second author was the lead author of the institutional grant from the Irvine Foundation that started the program.} Every summer five faculty members (mostly from within Pomona's Department of Mathematics \& Statistics) teach math in the program.  While the PAYS program is a college-wide program, it has {\bf directly} and positively affected the climate in the math and stats department in many ways, including the following.
 
 \begin{itemize}
 \item\ The myriad faculty who, over the years, have taught in the program have had a chance to think deeply about questions of inclusion and access as well as the value of community and collaboration. They have also gained a better appreciation of the lived experiences of our students of color, low-income students, and first generation students and the grit and tenacity that they bring to the classroom.
 \item\ Many of our math students, including our minoritized math students, participate in the PAYS program as TAs and mentors for the high school participants. The college students work closely with faculty and are trained as facilitators of collaborative learning. For the college students, the experience is worthwhile and socially meaningful. Moreover, the work is related to, and, in fact, requires, their mathematical training. These students feel connected to the department and its mission, and they are often departmental student leaders. 
 \item\ Regularly, PAYS students matriculate at Pomona College, and some have become math majors. Not only has  the presence of these students helped diversify the department, they have also brought the collaborative and inclusive culture of PAYS to the department.
 \end{itemize}
 
Through PSM, LiCMath, and PAYS, we endeavor to create a safe space for students of color, low-income students, and first generation students where their experience, their voice, and their story takes center stage.  At the same time, we want to involve as broad a swath of our students as possible for a number of reasons. First, without some work on the part of the (typically white) students who already feel safe in mathematics, it is hard to imagine a positive change in the climate. Second, we are concerned about stigmatizing a single group of students further, and third, rightly or wrongly, our assessment was that an approach that used a broader language would give us more room to operate. As a result, we have not been entirely consistent about ``who" the programs are for. For PSM, our targeted students are students of color, first-generation students, and low income students. In the PSM space, we often explicitly talk about race, gender, and ethnicity, and we try to center the experiences of the students involved.  In contrast, for learning communities, while we specifically encourage students of color, first-generation students, and low income students to participate, we invite all students to join.

\section{Departmental Involvement}

Like any other department, over the years, we have had disagreements over priorities. Even so, one of our shared principles has been to empower individual faculty to pursue their passions. As an example, we don't choose a common textbook or agree on a common syllabus for our calculus courses. The hope is that what we gain in faculty enthusiasm outweighs the advantages of common standards. At no time has any department member been directly obligated to participate in diversity initiatives, and, as you may expect, everyone participates at different levels. Even so, at this juncture, the whole department is supportive of the inclusion efforts, and, as a result, in addition to the programs outlined above, our faculty contribute to the sense of community and to our diversity work in myriad of other ways. For example, many routinely incorporate meals or social time with their students or with their research groups, and in choosing students for research and for leadership positions issues of openness and access are always considered.

The consensus over the importance of issues of access is the result of discussions---sometimes very heated---over a long period of time. Our first point of agreement---many years ago and with some pressure from our administration---was to agree that issues of access mattered in our job searches for new faculty. We have been including some version of the following in our job ads and in our discussions of individual candidates:
\begin{quote}
The department has directed much effort in creating a supporting community for all students, and is particularly interested in candidates who have experience working with students from diverse backgrounds and a demonstrated commitment to improving access to and success in higher education for students of color, low-income students, and first generation students. 
\end{quote}
Our deliberate discussion of diversity and access in the hiring process has meant that most of our current colleagues were hired with a clear understanding that this work is important to the department. As a \emph{PRIMUS} Associate Editor pointed out to us in the refereeing process, thinking of Diversity, Equity, and Inclusion work as a voluntary add-on to our jobs is actually one of the sources of the problems we face.  While we don't have a mechanism for distributing the work, and we don't force anyone to participate in any of the specific programs outlined here, most everyone in the department does see the DEI work as an important part of what our department does, and many contribute to it in one way or another.

\begin{figure}[H]
	\begin{center}
		\includegraphics[scale=.5]{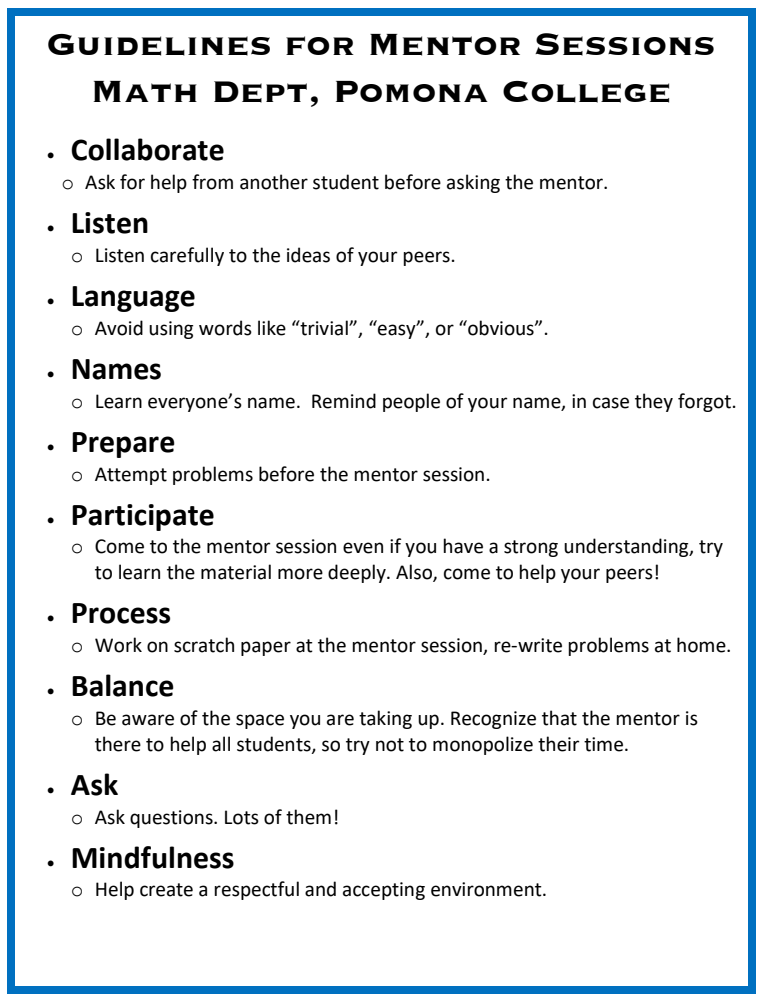}
\end{center}
	\caption{Mentor Session Guidelines -- posted in every classroom as a large poster. The original set of guidelines was created by one of our student mentors, Chris Donnay class of 2018.  \label{fig:mentor}}
\end{figure}

Among many other department initiatives, we also want to highlight three other efforts by our colleagues. The first is a community seminar on identity, culture, and mathematics organized in the Spring of 2018 by our colleague Gizem Karaali. The weekly seminar---with no credits or grades and open to all---brought together a group of students and faculty from our department. Each week the seminar discussed a new issue of interest to the community based on selected readings. Topics included \emph{Math Anxiety}, \emph{Microagressions}, \emph{Mathematics and Class}, \emph{Mathematical Inqueery}, \emph{Critical Mathematics Education}, \emph{Intentional Learning Communities and cohorts}, \emph{Role Models in Math}, \emph{Algorithmic Bias}, and \emph{Chinese Math, why do we think proofs are Greek?} In addition to providing a venue for in-depth discussion between faculty and students, the seminar came up with a number of action items. For example, the students, with faculty input, came up with a list of guidelines for our mentor sessions (see Figure \ref{fig:mentor}). The guidelines are posted in every classroom as large posters.

The second is the summer undergraduate research program PRiME\footnote{https://ehgoins.wixsite.com/prime} run by our colleague Edray Goins. Combining funding from the NSF and the College, Edray's program targets groups traditionally underrepresented in algebraic geometry--such as women and minoritized students, and, every summer, brings to campus a number of talented undergraduates who, along with a number of Pomona College students, do mathematical research. While it is too early to herald the results and impact on our department, it has been widely documented that undergraduate research is one of the best practices for retaining students in mathematics, and having that opportunity locally is invaluable.

The third is the EDGE program for women in mathematics\footnote{https://www.edgeforwomen.org/}, co-directed by our colleague Ami Radunskaya and hosted in Claremont {\bf three} times in recent years. While this program is for students who are already heading to graduate school, its availability, the programming around it, and the energy and enthusiasm of the participants is not lost to our undergraduate students.

While the connection between the myriad initiatives may not be immediately clear, we believe that the synergy across the programs is vital to the success of all of them.  Our students see the faculty involved and committed to diversity, equity, and inclusion.  Indeed, some of our students become involved in the satellite initiatives (even after graduating).  Our efforts do not (and cannot) have a directly linear trajectory or a simple cause and effect path, but we recognize and point out that our individual efforts have effectively combined to create many different spaces and venues where our students feel supported.

\section{Results}

Unquestionably, the Department of Mathematics \& Statistics at Pomona College looks and feels different now than it did when we started teaching (in 2002 and 1989, respectively).  The ongoing discussions we have throughout the building:  in department meetings, in our classrooms, in the hallways, and one-on-one with students tell us that, while we have not found a magic bullet, we are headed in the right direction.  We had zero Black math majors from 2003 to 2011 (when there were roughly 19 math majors per year); in 2020 out of our 45 graduating math majors, 5 (11.1\%) were Black and 5 (11.1\%) were Hispanic.  It is also worth noting that in the same class of 2020 math majors, 19 of 45 were women (42.2\%).  During a time in which the math major has been growing substantially, we are happy that the number of women majors has increased correspondingly while the total percentage of students of color has increased noticeably.

We are fortunate to work at an institution that is committed to higher educational access.  Over the years since initiating our programs, the overall student body at the College has diversified quite a bit.  We have seen parallel trends in the math and stats department.  The sample sizes are too small to do full comparative analyses, but we are encouraged by the direction of the trend, particularly with students of color becoming math majors.  The institution has long been approximately half women.  Unfortunately, in the last two decades, we have seen  little to no change in the proportion of women mathematics majors, suggesting that we should possibly think even more carefully about our female students and likely our female students of color.  We hope that our initiatives are part of the mechanism encouraging women and marginalized students to pursue math.  However, without a better experimental design, we are not in a position to claim that our programs are causing the changing demographics in the math and stats department.

\subsection{Qualitative Student Feedback \label{kayla}}

As they were graduating, we surveyed the inaugural PSM class (c/o 2018).  Here, we provide their qualitative feedback as a way to give a sense of their experiences with PSM.

\begin{itemize}
\item\ ``Having PSM made it easier to acknowledge when I was struggling the most."
\item\ ``This group is probably the main reason I chose math."
\item\ ``Thank you so, so much.  This community really keeps me going.  Isolation is easy to come by at Pomona."
\item\  ``More open/safe conversations...  re: imposter syndrome/stereotype threat."
\item\ ``I had a particularly rough time my first semester because only 2 PSM members were in my calc class and they both dropped, leaving me alone.  I know it's difficult to ensure there will be PSM members together since you cannot force people to take a course, but maybe some form of additional communication could better help facilitate that."
\item\ ``[PSM] is very math-centric which mades sense b/c the word math is in there but it might push some away from being in it."
\item\ ``the relationships that I formed through this group as both a mentor and as a mentee have allowed me to be successful as a low-income first-gen woman in STEM."
\item\ ``PSM could improve by discussing non-academic topics during meetings and one-on-one advising.  Many of the PSM peers, myself included, are the first in their family to attend college.  This whole experience, the different social rules, and the full new atmosphere are brand new.  I think that spending more time discussing this and potential responses / coping strategies / lessons learned would be helpful."
\item\ ``PSM was great and I really really appreciated it.  It truly supported me and helped me during my first year."
\end{itemize}

The quantitative information we have on diversification of our department is neither comprehensive nor striking enough to make bold claims about having figured out the perfect blend of programs to create an ideal space for all of our students.  However, the qualitative information and our sense for how the department has changed over the last two decades {\bf does} give us hope that we are moving in the right direction.  Our students seem to be well served by the programs, and we can tell that their engagement is likely due in a large part to our engagement.  That said, we are well aware that it is imperative for us to reach each new group of students as they matriculate at Pomona College.  Even after the programs have been set up and put into motion, we are obligated to continually revisit and revise their effectiveness.

We have stayed in touch with many of the students who have passed through the programs described above.  Some of the students are in graduate school or finished with PhDs or master's degrees. Others are working in tech or finance or education.  They are a heterogeneous group  who have taken their Pomona educations in many different directions.  In our efforts to bring their voices into this manuscript, we have reached out to our alumni and reconfirmed the variety of student experiences.   We certainly continue to hear from some of our alumni who are still navigating their own identities in mathematics, and we sometimes hear from students that their voices and their experiences are not sufficiently valued in the math and stats department and that they continuously feel that they have to carry a much heavier burden than their peers.  Other students, however, are able to directly connect the efforts from our programming to their post-Pomona trajectory as evidenced by the recent reflection from a PSM student who is currently in a PhD program:

\begin{quote}
I think Pomona math's diversity programming made a huge difference in the way that I conceived of and conceptualized my mathematical self at the time. I gained so much confidence, I was comfortable being wrong and asking for help, and I trusted myself to contribute to collaborative settings.  This part of me wilted and at one point had disappeared at [University XXX]\footnote{A top research institution that has been redacted for privacy reasons.}. I really struggled to feel competent or even smart during the first two years of my PhD, when none of this support was available to me, in a program with $<$15\% women and $<$5\% POC (not including east-Asian backgrounds), being constantly undermined by peers and professors, interrupted, ignored in group settings. I really leaned on the spirit of Pomona's math community from my past to figure out how to handle this.  ...[T]he foundation was available for me to have just enough disbelief in the new voices [University XXX] had cultivated in me (I am stupid, I don't belong here, I faked my way in, no one takes me seriously) to talk back to them and wrangle them out of the driver's seat.
 \end{quote}

\section{Challenges}

In recent years, we have put quite a bit of effort into issues of diversity and access. Our department indeed feels and looks different than earlier times, and we are excited about the changes. Even so, many challenges remain. Here are a scattering of questions with which we continue to wrestle.

\begin{itemize}
\item What are best practices in targeting students for our programs, and how do we avoid unintentionally siloing students and creating new barriers between different groups of students?
\item\ How much of a priority should diversity issues be? Are there aspects of the work that should be considered a core part of our mission? If so, how should the department support these efforts?
\item\ How do we meaningfully assess what we do? In the past several years, we have thrown everything and the kitchen sink at the lack of diversity among our students.  We have seen some improvement, but are the changes due to our efforts or to the changes in demographics brought about by the Pomona College admissions staff? Which aspects of our programs are most effective? Some of what we do is incredibly labor intensive. Do the outcomes justify the energy spent?
\item\ How do we put a structure in place which assures that we learn from past experience, that we pass along best practices, and that we achieve our goals efficiently? As an example, what kind of training makes sense for our faculty and student mentors?
\item\ How can we learn from the vast literature on education research? We do not have the in-house expertise to effectively evaluate what we do and to separate out the essential from the superfluous.
\end{itemize}

\section{Takeaways}

As we said from the onset, neither one of us are experts in research on education or in the cutting edge debates on diversity in the academy. We are run of the mill mathematicians and statisticians who believe that we have to play a role. In this decades long process, we ourselves have been learning and often have adjusted our thinking, our language, and our approaches as we gained more understanding. Our main take aways from this experience are: 
\begin{enumerate}
\item
There is no sliver bullet, 
\begin{quote} try things!  try implementing our ideas.  try implementing your ideas.  ask your students what they think would work for them.\end{quote}
\item
change takes time, 
\begin{quote} time is likely the most important factor.  not only is change slow, but it also requires sustained commitment.  that is, we spend many hours on the programming described above.  it is difficult to outsource or scale the work.  you just need to put in the time. \end{quote}
\item
it is possible and necessary to make progress, 
\begin{quote} anything is better than nothing.  whatever you do will be appreciated, and the effects of your programs will infuse into all aspects of your own work and community.  \end{quote}
\item
everyone plays a role, 
\begin{quote} the work isn't only done by the faculty members.  one of our big keys to success was our long time administrative assistant who was able to help create community in myriad ways.  another group who led to our success were older or graduated students.  we were able to lean on their experiences as examples for our more junior students.\end{quote}
\item
but not everyone needs to be on board to start, and 
\begin{quote} it's okay (and expected!) if your department has different levels of commitment to this work.  join with the people in your department who want to do the work, and you may find your group expanding over time. \end{quote}
\item
no progress is permanent, you have to keep reproducing it.
\begin{quote} we have a long way to go, and we are in it for the long haul.  we hope that you are, too. \end{quote}
\end{enumerate}

In conclusion, there is no question that our efforts to create an inclusive community in the Department of Mathematics \& Statistics will be and should be an ongoing pursuit. We must continue to adapt our programs and our approach to the changes in our student body and to the challenges and uncertainties posed by the outside world. Our students have been a source of inspiration, and we are impelled to engage them and to learn from them. We are also happy to share our experiences in greater detail. Feel free to contact us, especially if you are interested in starting a program and have logistical questions.

\newpage

\bibliographystyle{plainnat}
\bibliography{CCCbib}

\end{document}